\documentclass[12pt, reqno, psamsfonts]{amsart}

\usepackage{latexsym, amsfonts, amsmath, amsthm}
\usepackage[psamsfonts]{amssymb}

\newtheorem{theorem}{Theorem}[section]
\newtheorem{lemma}[theorem]{Lemma}

\newtheorem{observation}[theorem]{Observation}
\newtheorem{proposition}[theorem]{Proposition}
\newtheorem{corollary}[theorem]{Corollary}

\theoremstyle{definition}
\newtheorem{definition}[theorem]{Definition}
\newtheorem{example}[theorem]{Example}

\theoremstyle{remark}
\newtheorem{remark}[theorem]{Remark}

\numberwithin{equation}{section}

\newcommand{\pair}[2]{(#1|#2)}  \newcommand{\paircd}{(\cdot|\cdot)}

\newcommand{\bb}{{\mathbf b}}

\newcommand{\bh}{{\mathbf h}}

\newcommand{\ICO}[1][\CO]{{\mathfrak I(#1)}}

 \newcommand{\del}{\partial}

\newcommand{\mapright}[1]{\smash{\mathop{\longrightarrow}\limits^{#1}}}
\newcommand{\mapleft}[1]{\;\smash{\mathop
   {\longleftarrow}\limits^{#1}}}

\newcommand{\iso}{\;\smash{\mathop{\longrightarrow}\limits^{\sim}}\;}

\newcommand{\inv}{^{-1}}
\newcommand{\ovl}{\overline}
\newcommand{\wt}{\widetilde}

\newcommand{\half}[1][1]{\frac{#1}{2}}
          \newcommand{\g}{{\mathfrak g}}

    \newcommand{\kk}{{\mathfrak k}}

          \newcommand{\p}{{\mathfrak p}}
              \newcommand{\s}{{\mathfrak s}}

\newcommand{\fsp}{\mathfrak{sp}}
\newcommand{\fsl}{\mathfrak{sl}}

\newcommand{\Z}{{\mathbb Z}}
\newcommand{\N}{{\mathbb N}}    \newcommand{\hN}{{\half\mathbb N}}

\newcommand{\C}{{\mathbb C}}

  \newcommand{\Om}{\Omega}

\newcommand{\ga}{\gamma}

\newcommand{\la}{\lambda} 

\newcommand{\sig}{\sigma}  
\newcommand{\vsig}{\varsigma}
\newcommand{\vth}{\vartheta}

\newcommand{\B}{{\mathcal B}}
\newcommand{\D}{{\mathcal D}}

\newcommand{\JJ}{{\mathcal J}}

                                          \newcommand{\OO}{{\mathcal O}}

                                         \newcommand{\T}{{\mathcal T}}
                                           \newcommand{\U}{{\mathcal U}}
\newcommand{\V}{{\mathcal V}}
                                        \newcommand{\W}{{\mathcal W}}
\newcommand{\X}{{\mathcal X}}
\newcommand{\Y}{{\mathcal Y}}

\newcommand{\gr}{\mathop{\mathrm{gr}}\nolimits}

\newcommand{\End}{\mathop{\mathrm{End}}\nolimits}
\newcommand{\Hom}{\mathop{\mathrm{Hom}}\nolimits}

\newcommand{\tr}{\mathop{\mathrm{tr}}\nolimits}

\newcommand{\gog}{\g\oplus\g}
\newcommand{\sos}{\s\oplus\s}

\newcommand{\Ug}[1][]{{\mathcal U}_{#1}({\mathfrak g})}

\newcommand{\Sg}[1][]{S^{#1}({\mathfrak g})}

\newcommand{\Kill}[1]{({#1})}

\newcommand{\SpL}{Sp(L,\C)}  \newcommand{\spL}{\fsp(L,\C)}

\newcommand{\CO}{\ovl{\OO}}
\newcommand{\GtS}{G\times S}  
\newcommand{\INV}{^{inv}}

\newcommand{\RR}[1][]{{\mathcal R^{#1}}}
\newcommand{\PP}[1][]{{\mathcal P^{#1}}}
\newcommand{\II}[1][]{{\mathcal I^{#1}}}

\newcommand{\PV}{\PP\INV}  \newcommand{\IV}{\II\INV}
\newcommand{\MV}{\M\INV}
\newcommand{\WV}{\W\INV}
\newcommand{\We}{\W^{even}}    \newcommand{\Wo}{\W^{odd}}
\newcommand{\Pe}{\PP^{even}}   \newcommand{\Po}{\PP^{odd}}
\newcommand{\Ie}{\II^{even}}

\newcommand{\taug}{\tau_{\g}}
\newcommand{\tauB}{\tau_{\B}}

\newcommand{\vthg}{\vth_{\g}}
\newcommand{\vthB}{\vth_{\B}}

\newcommand{\xiB}{\xi_{\B}}

\newcommand{\M}{\mathcal M}

\newcommand{\Ct}[1][]{{\C[t]}}

\newcommand{\PIV}{\PP\INV/\II\INV}

\newcommand{\Spc}{Sp(L)}     \newcommand{\spc}{\fsp(L)}
\newcommand{\Sc}{S_{c}}     
\newcommand{\Gc}{G_{c}}    \newcommand{\gc}{\g_{c}}

\newcommand{\red}{^{red}}

\newcommand{\PBO}{\{\cdot,\cdot\}}

\newcommand{\PB}[2]{\{#1,#2\}}

\newcommand{\txi}{\wt{\xi}}

\newcommand{\gsh}{\g^{\sharp}}

\newcommand{\oY}{\ovl{\Y}}

\newcommand{\Ei}{E_{\infty}} \newcommand{\Zi}{Z_{\infty}}
\newcommand{\Bi}{B_{\infty}}

\newcommand{\leftsuper}[1]{{}^{#1}{\kern -2pt}}

\begin{document}

\title
{Dixmier Algebras  for  Classical Complex Nilpotent
Orbits via Kraft-Procesi Models I}
\author{Ranee Brylinski}
\address{Department of Mathematics,
        Penn State University, University Park 16802}
        \address{current mailing address (2002-...):
        P.O.Box 1089, Truro MA 02666-1089}
        \address{}
\email{ rkb248@yahoo.com} \urladdr{www.math.psu.edu/rkb}
\thanks{}

\dedicatory{Dedicated to Professor Alexander Kirillov on his $2^6$ birthday}

\begin{abstract}
We attach a Dixmier algebra $\B$ to the closure $\CO$ of any nilpotent orbit of
$G$ where $G$ is $GL(n,\C)$, $O(n,\C)$ or $Sp(2n,\C)$.
This algebra $\B$ is a noncommutative analog of the coordinate ring
$\RR$ of  $\CO$, in the sense that $\B$ has a $G$-invariant algebra filtration
and $\gr\B=\RR$.

We obtain $\B$ by making a noncommutative analog of
the Kraft-Procesi construction which modeled  $\CO$ as the  algebraic symplectic
reduction of a finite-dimensional symplectic vector space $L$.
Indeed $\B$ is a subquotient of the Weyl algebra for $L$.

$\B$ identifies with the quotient of $\Ug$ by a two-sided ideal $J$, where $\g=Lie(G)$.
Then $\gr J$ is the ideal $\ICO$ in $\Sg$ of  functions vanishing on $\CO$.
In every case where $\OO$ is connected, $J$ is a completely prime primitive ideal.
\end{abstract}

\maketitle

\section{Introduction}
\label{sec:intro}

By means of symplectic reduction in the setting of complex algebraic varieties,
Kraft and Procesi (\cite{KP1},\cite{KP2})  constructed a model  of
the closure   of  any  nilpotent coadjoint orbit $\OO$  of $G$
when $G$ is one of  the classical groups
$GL(n,\C)$, $O(n,\C)$ and $Sp(2n,\C)$.
The symplectic aspect is   not actually mentioned, but
the construction is clearly symplectic.

In this paper we give a noncommutative analog, or quantization,
of the Kraft-Procesi  construction.
The result is that we attach a Dixmier algebra $\B$ to each orbit closure $\CO$.
Our algebra $\B$ has a $G$-invariant   algebra filtration  and we show that
$\gr\B$ is isomorphic, as a graded Poisson algebra, to the coordinate ring
$\RR$ of $\CO$.

In fact, $\B$ identifies, as a filtered algebra, with the   quotient $\Ug/J$ of the
universal enveloping algebra $\Ug$ of $\g=Lie(G)$ by some  two-sided ideal $J$.
Then $\gr J$ is the
ideal $\ICO$ in $\Sg$ defining $\CO$. We find  that $J$  is stable under the
principal anti-automorphism of $\Ug$, and also under the anti-linear automorphism
of $\Ug$ defined by a  Cartan involution of $\g$.

The Kraft-Procesi  construction attaches to $\OO$ a complex symplectic vector
space $L$ together with a Hamiltonian action of $\GtS$ on $L$, where $S$
is an auxillary complex reductive Lie group. The actions of $G$ and $S$ lie inside
the symplectic group $\SpL$. Kraft and Procesi
show   that $\CO$ is scheme-theoretically the algebraic
symplectic reduction of $L$ by $S$.
In this way, they obtain $\RR$  as a subquotient of the algebra $\PP$ of polynomial
functions on $L$. More precisely, $\RR$ is realized as $\PIV$ where
$\II$ is an ideal in $\PP$ and the superscript $inv$ denotes taking $S$-invariants.

It is nice from the viewpoint of  representation theory to regard $\RR$
as a subquotient of the algebra  $\Pe$ of even polynomials. (We can do this
since $\PV$ lies in $\Pe$.)  For $\Pe$ is the coordinate ring of the closure
of the minimal nilpotent orbit $\Y$ of $\SpL$.

To  make a   noncommutative analog  of the Kraft-Procesi construction,
we start from the fact that  there is a unique Dixmier algebra attached
to $\oY$, namely the quotient of $\U(\spL)$ by its Joseph ideal $\JJ$.
We can model $\U(\spL)/\JJ$ as the even part $\We$ of the Weyl algebra
$\W$ for $L$ and then $\gr\We=\Pe$.    There is an obvious
quantization of the Hamiltonian $S$-action on $L$, namely the  natural
$(\sos,\Sc)$-module structure on $\W$.
Here the subscript $c$ denotes taking a compact real form.

We define $\B$ to be the coinvariants  for the $(\sos,\Sc)$-action on $\We$.
A priori, $\B$ is a $(\gog,\Gc)$-module with a $\Gc$-invariant filtration,
but $\B$ is not an algebra. However, we easily identify $\B$ as the quotient
by a two-sided ideal of $\WV$ (where the superscript again indicates
taking $S$-invariants, or  equivalently, $\Sc$-invariants).
In this way, $\B$ becomes a  filtered algebra and a subquotient   of $\We$.

Our main result (Theorem \ref{thm:B})  is to compute the associated graded algebra
$\gr\B$.  It is easy to see that $\gr\B$ is some quotient of $\PIV$, but in
fact we prove $\gr\B=\PIV$.  To do this, we recognize $\B$ as the degree zero part
of the relative Lie algebra homology  $H(\sos,\Sc;\We)$.
We consider the   standard complex which computes this homology,
introduce a filtration and then apply   the spectral sequence for a filtered complex.
We compute the $E_1$ term of the spectral sequence by using
the fact proven by Kraft and Procesi that $\II$ is a complete intersection ideal.
Then we easily show   $E_1=\Ei$.

We establish some properties of $\B$ and the corresponding ideal $J$.
If $\OO$ is connected then $J$ is a completely prime primitive ideal
(Corollary \ref{cor:B_J}).
In every case, $\B$ admits a unique
$\gsh$-invariant Hermitian inner product $\paircd$  such that $\pair{1}{1}=1$,
where $\gsh$ is a real form of $\gog$ with $\gsh\simeq\g$
(see Proposition \ref{prop:3.1}).
This prompts the question as to whether
$J$ is  ``good" in the sense that  $J$ is maximal and $\Ug/J$ is unitarizable.
The  latter property means (since $\B\simeq\Ug/J$)  that $\paircd$ is  positive-definite.

Attaching ``good" ideals to $\OO$ is an important problem in representation
theory and the orbit method.   Quite a bit of work has been done on this (see
e.g. some of the references and authors cited below) but the  problem
for nilpotent orbits of a complex semisimple Lie group remains unsolved.

If $G=GL(n,\C)$, then our $J$ is  good  (see    Remark  \ref{rem:GL}(i) and \cite{Brylinski:Dix_KP_II} ).
But $G=GL(n,\C)$ is really a very special case for us as the geometry of $\OO$
is incredibly nice, including but not limited to the fact that $\CO$ is always normal.
For   $G=O(n,\C)$ or $G=Sp(2n,\C)$, it is not the case that $J$ is always good.
Certainly if $\CO$ is not normal, we should not expect $J$ to be good.

Our point of view (to be justified in \cite{Brylinski:GQ_KP}) is that $\B$ is the
``canonical" quantization of the Kraft-Procesi construction, and so the failure
of $J$ to be good is really a statement about  $\CO$. The next step is then to
investigate whether we can make a modification to our quantization process
in order to obtain some good ideals in $\Ug$ attached to $\OO$ (and even
its covers).

This paper is the first in a series. In the subsequent papers we
make explicit the important role of  Howe duality in our project.
Indeed, $\We$ is the Harish-Chandra module of the (even)
oscillator representation of $Sp(L,\C)$, and the pair $(G,S)$
constitutes a sequence of Howe dual pairs (see Remark
\ref{rem:Howe}). In taking coinvariants, we are implementing a
sequence of  Howe  duality ``operations".  Each ``operation"  is
like implementing a Howe duality correspondence,  except  that we
do \emph{not} pass to the the irreducible quotient.     In
working on this project (which started in earnest in  the summer
of 2001  -- and is part of a program we began in 1994), we have
been reading the Howe duality literature. We have been influenced
by especially  the papers \cite{Howe}, \cite{Moeglin_howe},
\cite{AdamsBarbasch} and  \cite{Li}.

Our first construction of
the Dixmier algebra $\B$ actually came out of the ideas of Howe duality
and  quantization by constraints.  This is given in \cite{Brylinski:GQ_KP} and lies
more in the realm of harmonic analysis than algebra.
Our starting point there  is the fact (\cite{KS}) that
$L$ is hyperk\"ahler and the Kraft-Procesi construction is the algebraic analog
of the  hyperk\"ahler reduction of $L$  by $\Sc$.

The notion of Dixmier algebra  for nilpotent orbits
(including their  closures and their coverings)
was first developed in work  of McGovern, Joseph and Vogan.
See e.g.  \cite{McG:mem}, \cite{Vog1},  and \cite{Vog2}.
The motivation for these authors and for most  Dixmier algebra
theorists is the search for completely prime primitive ideals.
This motivation is very important for us too;  we also find
additional  motivations  coming from  star products and from
geometric quantization.

The results  in this paper should be compared with the work in
\cite{AdamsBarbasch}, \cite{Barbasch}, \cite{Barbasch2},
\cite{BarbaschV1}, \cite{BarbaschV2}, \cite{Brylinski:EDQ},
\cite{CG}, \cite{DLO}, \cite{Duflo}, \cite{Gross}, \cite{LS},
\cite{Li}, \cite{McG:yale}, \cite{McG:mem}, \cite{Moeglin_gln},
\cite{Moeglin_howe} and \cite{Vog3}.  Some of this comparison work
will be done in \cite{Brylinski:Dix_KP_II}  and
\cite{Brylinski:GQ_KP}.

Part of this work was carried out while I was visiting the
IML and the CPT of the Universit\'e de la M\'editerran\'ee
in the summer of 2001, and I thank  my colleagues
there for their hospitality. I especially thank Christian Duval
and Valentin Ovsienko for some very  valuable discussions.

It is a real pleasure to dedicate this article to Sasha Kirillov
whose discoveries have  opened up so many new vistas, starting of course
with the Orbit Method. I warmly thank  him for his friendship and his
interest  in my own work.

\section{Dixmier Algebra for the closure  of a complex nilpotent  orbit}
\label{sec:OO}

Let $G$ be a reductive complex  algebraic group.
Let  $\g$ be the Lie algebra of $G$ and let $\g^*$ be the dual space.
Then $G$ acts on $\g$ and $\g^*$ by, respectively,
the adjoint action and  the coadjoint action.
The symmetric algebra $\Sg=\oplus_{p=0}^\infty \Sg[p]$
is the algebra of polynomial functions on $\g^*$.
The $G$-invariants  form the graded subalgebra
$\Sg^G=\oplus_{p=0}^\infty \Sg[p]^G$.
We can fix some  nondegenerate $G$-invariant bilinear form
$\Kill{\cdot,\cdot}$ on $\g$.

The \emph{nullcone} in $\g^*$ is the set of $\la\in\g^*$ which satisfy the
following equivalent properties:

(i) the closure  of the coadjoint orbit of $\la$ contains zero

(ii) the coadjoint orbit of $\la$ is stable under dilations of the vector space $\g^*$

(iii) every nonconstant homogeneous $G$-invariant  in $\Sg$ vanishes on $\la$

(iv) $\la=\Kill{x,\cdot}$ where $x$ is a nilpotent in $\g$

\noindent The nullcone is $G$-stable and breaks into finitely many   orbits
of $G$, which are then called the \emph{nilpotent coadjoint orbits}, or
simply  the \emph{nilpotent orbits},  of $G$.

Let $\OO$ be a nilpotent  orbit of $G$.
The closure $\CO$  is a  complex algebraic subvariety of  $\g^*$;
but $\CO$ may be reducible if $G$ is disconnected.
The coordinate ring $\C[\CO]$ of $\CO$ is the quotient algebra
\begin{equation}\label{eq:R=}
\RR=\Sg/\ICO
\end{equation}
where $\ICO$ is  the ideal  of functions which vanish on $\CO$.
Then $\ICO$ is a graded ideal and
$\RR=\oplus_{p=0}^\infty \RR[p]$ is a graded algebra where
$\RR[p]=\Sg[p]/\ICO^p$.
Each space $\RR[p]$ is a finite dimensional completely reducible
representation of $G$.
Kostant's description of $\Sg$ as a module over  $\Sg^G$
implies that  all $G$-multiplicities in $\RR$ are finite.

$\RR$ inherits from $\Sg$ the structure of a graded Poisson algebra where
$\PB{\RR[p]}{\RR[q]}\subseteq\RR[p+q-1]$.
This  Poisson bracket $\PBO$ on $\RR$ is $G$-invariant and corresponds
to  the  holomorphic Kirillov-Kostant-Souriau  symplectic form on $\OO$.

In this situation,  we define Dixmier algebras in the following way.
We fix a   Cartan involution $\vsig$ of  $\g$. Then  $\vsig$
corresponds to a compact real form $\Gc$ of $G$ with Lie algebra $\gc$.
Let $\N=\{0,1,2,\dots\}$.
\begin{definition}\label{def:Dix}
A \emph{Dixmier algebra} for $\CO$ is a quadruple $(\D,\xi,\tau,\vth)$ where
\begin{itemize}
\item  $\D$ is a filtered algebra with an increasing algebra filtration
$\D=\cup_{p\in\N}\,\D_p$ such that $\gr\D$ is commutative.
\item   $\xi:\g\to\D_1$, $x\mapsto\xi^x$, is a  homomorphism of Lie algebras
and $\xi$ induces an isomorphism  of graded Poisson algebras from
$\Sg/\ICO$ onto $\gr\D$.
\item  $\tau$ is a  filtered algebra anti-involution of $\D$ such that
$\tau(\xi^x)=-\xi^x$.
\item  $\vth$ is an anti-linear filtered algebra involution
such that  $\vth(\xi^x)=\xi^{\vsig(x)}$.
\end{itemize}
\end{definition}
Here are some explanations  about the definition.
First, commutativity of $\gr\D$ implies that $\gr\D$ has a natural structure of
graded  Poisson algebra; here  the commutator in $\D$ induces  the
Poisson bracket on $\gr\D$.
Second,  $\xi$ extends naturally to a filtered algebra homomorphism
\begin{equation}\label{eq:hxi}
\txi:\Ug\to\D
\end{equation}
Let $J$ be the kernel of $\txi$. Then $\gr J$ is a Poisson ideal of $\Sg$,
and $\gr\txi$ induces a $1$-to-$1$ homomorphism $\zeta:\Sg/\gr J\to\gr D$
of graded Poisson algebras.  We require  that
$\zeta$  is surjective and
\begin{equation}\label{eq:grJ=}
\gr J=\ICO
\end{equation}
Notice that $\zeta$ is surjective  if and only if  $\txi$ is surjective in each
filtration degree; then $\txi$  induces a filtered algebra  isomorphism
\begin{equation}\label{eq:isoD}
\Ug/J\iso\D
\end{equation}
Third,  $\tau$ satisfies $\tau(c A)=c\,\tau(A)$, $\tau(A+B)=\tau(A)+\tau(B)$,
and $\tau(AB)=\tau(B)\tau(A)$ where $A,B\in\D$ and $c\in\C$.
Fourth, $\vth$ satisfies $\vth(c A)=\ovl{c}\,\vth(A)$, $\vth(A+B)=\vth(A)+\vth(B)$,
and $\vth(AB)=\vth(A)\vth(B)$. Clearly $\tau\vth=\vth\tau$.

Notice  that $\D$ and $\xi$ (and $\vsig$) uniquely determine
$\tau$ and $\vth$,  if the latter exist.
 Indeed,  the endomorphisms $x\mapsto-x$ and $x\mapsto\vsig(x)$ of $\g$
extend uniquely  to $\taug$ and $\vthg$, where $\taug$ is  an algebra
anti-involution  of $\Ug$ and $\vthg$ is an antilinear algebra involution of $\Ug$.
Then, via (\ref{eq:isoD}), $\taug$ and $\vthg$ induce $\tau$ and $\vth$.

We have  an obvious notion of isomorphism of Dixmier algebras for $\CO$:
$(\D,\xi,\tau,\vth)$ is isomorphic to $(\D',\xi',\tau',\vth')$ if there is a
filtered algebra isomorphism $\eta:\D\to\D'$ such that
$\xi'=\eta\circ\xi$,  $\eta\circ\tau=\tau'\circ\eta$, and
$\eta\circ\vth=\vth'\circ\eta$.
We can easily classify Dixmier algebras.
\begin{observation}\label{obs:Dix}
Suppose $(\D,\xi,\tau,\vth)$ is a Dixmier algebra for $\CO$.
In addition to \textup{(\ref{eq:grJ=})}, $J$ satisfies
\begin{equation}\label{eq:J..}
\taug(J)=J \qquad\textrm{and}\qquad \vthg(J)=J
\end{equation}
In this way, we get a bijection between
\textup{(}isomorphism classes of\textup{)} Dixmier algebras for $\CO$
and two-sided ideals $J$ of $\Ug$ satisfying
\textup{(\ref{eq:grJ=})} and \textup{(\ref{eq:J..})}.
\end{observation}
\begin{proof}
Clearly (\ref{eq:grJ=}) and (\ref{eq:J..}) imply that
$(\Ug/J,\iota,\taug',\vthg')$  is a Dixmier algebra for $\CO$,
where $\iota$ is the obvious map and
$\taug'$ and $\vthg'$ are  induced by  $\taug$ and $\vthg$.
Conversely, if $(\D,\xi,\tau,\vth)$ is    given, then
$(\Ug/J,\iota,\taug',\vthg')$  is isomorphic to it  via (\ref{eq:isoD}).
\end{proof}

\section{Properties of Dixmier algebras}
\label{sec:3}

The  hopes in constructing a Dixmier algebra are
(i)  $J$ will be a completely prime primitive ideal of $\Ug$, or even better,
a completely prime  maximal  ideal, and (ii) $\D$ will be unitarizable.
See   \cite[3.1]{Dixmier} for the  definitions of the terms in (i).

To understand (ii),  we observe that  Definition \ref{def:Dix}  makes $\D$
into a  $(\gog,\Gc)$-module.
Indeed,   the natural  $(\gog,\Gc)$-module structure   on
$\Ug/J$ transfers over to $\D$ via (\ref{eq:isoD}).  Then
$\gog$ acts on $\D$ though the representation
\begin{equation}\label{eq:Pi}
\Pi:\gog\to\End\D,\qquad (x,y)\mapsto\Pi^{x,y}
\end{equation}
where $\Pi^{x,y}(A)=\xi^xA-A\xi^y$. The action of $\Gc$ corresponds to
the subalgebra   $\{(x,x): x\in\gc\}$.

Next consider the subalgebra $\gsh=\{(x,\vsig(x)): x\in\g\}$ of $\gog$.
We say $\D$ is \emph{unitarizable} if $\D$ admits a $\gsh$-invariant
positive definite Hermitian inner product.
In this event,  by a theorem of Harish-Chandra,  the  operators
$\Pi^{x,\vsig(x)}$  correspond  to a unitary representation of $G$ on the
Hilbert space  completion  of $\D$.
This unitary representation is then a quantization of $\OO$
in the sense of geometric quantization, if we view $\OO$  as
a \emph{real} symplectic manifold. (If $\C[\CO]\neq\C[\OO]$,
then this might be just  a piece of a quantization of $\OO$.)

Notice that the following three properties are equivalent:
(i) $J$ is maximal, (ii) $\D$ is a simple ring, and (iii) the representation $\Pi$
is irreducible.

Our formalism gives some partial results pertaining to hopes (i) and (ii).
Notice that   $\D_0=\C$ by  (\ref{eq:isoD}).
\begin{proposition}\label{prop:3.1}
Suppose $(\D,\xi,\tau,\vth)$ is a Dixmier algebra for $\CO$ and
$J=\ker\txi$.   Then
\begin{itemize}
\item[(i)]  $J$ has an infinitesimal character.
\item[(ii)]   If $\CO$ is irreducible then $J$ is a completely prime
primitive ideal  in $\Ug$.
\item[(iii)]  There is a  a unique $\Gc$-invariant  projection $\T:\D\to\C$.
This map $\T$ is a trace,  i.e., $\T(AB)=\T(BA)$.
\item[(iv)]    $\D$ admits a  unique $\gsh$-invariant Hermitian inner product
$\paircd$  such that $\pair{1}{1}=1$, and it is given by
\begin{equation}\label{eq:pair}
\pair{A}{B}=\T(AB^{\vth}).
\end{equation}
where $B^{\vth}=\vth(B)$.
\end{itemize}
\end{proposition}
\begin{proof}
(i)  This means (for any proper two-sided ideal $J$) that   $J$ contains a
maximal ideal of the center  of $\Ug$.    This happens
if and only $\gr J$ contains $\Sg[+]^{G}=\oplus_{p>0}\Sg[p]^{G}$.
But    $\gr J=\ICO$ and  $\ICO\supset\Sg[+]^{G}$ since
$\CO$ lies in the nullcone.
(ii)  If $\CO$ is irreducible then $\ICO$ is a prime ideal in $\Sg$ and
so $J$ is a completely prime ideal in $\Ug$. This together with (i)
implies, by a result of Dixmier,  that  $J$ is primitive.
(iii)  Since $\ICO\supset\Sg[+]^{G}=\Sg[+]^{\Gc}$,  we have $\RR^{\Gc}=\C$
and so  $\D^{\Gc}=\C$.
Thus we get a unique $\Gc$-invariant  projection map $\T$.
Now $\Gc$-invariance implies  $\T([\xi^x,A])=0$ where $x\in\g$.
We  can write this as $\T(\xi^xA)=\T(A\xi^x)$. Iteration gives
$\T(\xi^{x_1}\cdots\xi^{x_k}A)=\T(A\xi^{x_1}\cdots\xi^{x_k})$.
This proves $\T(BA)=\T(AB)$ since the $\xi^x$ generate $\D$.
(This is the same proof as in \cite[Proposition 8.4]{Brylinski:EDQ}.)
(iv)  Suppose $\paircd$ is an  inner product with the desired properties.
Then $\pair{A}{1}=\T(A)$.  Now $\gsh$-invariance   means that the operators
$\Pi^{x,\vsig(x)}$ are skew-hermitian, or equivalently,
$\pair{\xi^xA}{B}=\pair{A}{B\xi^{\vsig(x)}}$.
So for $B=\xi^{x_1}\cdots\xi^{x_k}$,
we have $\pair{A}{B}=\pair{\xi^{\vsig(x_1)}\cdots\xi^{\vsig(x_k)}A}{1}
=\pair{B^\vth A}{1}=\T(B^\vth A)$. The result is now clear.
 \end{proof}

\begin{corollary}\label{cor:3.2}
$\D$ is unitarizable if and only if the pairing defined by \textup{(\ref{eq:pair})}
is positive definite.
If $\D$ is unitarizable, then $J$ is maximal.
\end{corollary}
\begin{proof}
Both statements follow from the uniqueness in Proposition \ref{prop:3.1}(iv).
\end{proof}

\begin{example}\label{exam:min}
Suppose $\OO$ is the minimal nilpotent orbit in $\g$ where $\g$ is simple
and $\g\neq\fsl(2,\C)$. This is a case where  $\C[\CO]=\C[\OO]$.
Then there is a unique Dixmier algebra for $\CO$.
This follows by Observation \ref{obs:Dix}  since there is exactly one
choice for $J$ satisfying (\ref{eq:grJ=}) and (\ref{eq:J..}).
Moreover, (i) $J$ is a  completely prime maximal  ideal of $\Ug$,
and (ii) $\Ug/J$ is unitarizable if $\g$ is classical.

Indeed, there  is a unique two-sided  ideal $J$ satisfying (\ref{eq:grJ=})
and $\taug(J)=J$; see \cite[proof of Proposition  3.1]{AB:starmin}.
Since the ideal $\ICO$ is preserved by the antilinear algebra involution
of $\Sg$ defined  by $\vsig$,  it follows by
the uniqueness of $J$ that  $\vthg(J)=J$. Since $\CO$ is irreducible,
Proposition \ref{prop:3.1}(ii) implies that $J$ is completely prime.
If $\g\neq\fsl(n,\C)$, then $J$ is the Joseph ideal  and this is   maximal
by \cite[Theorem 7.4]{Jos}.  If $\g=\fsl(n,\C)$ ($n\ge 3$), then $J$ is maximal by
\cite{VdB}. Finally,   unitarizabilty is known; see \cite[Theorem 9.1]{AB:starmin}
for   a  uniform construction of these unitary representations on spaces of
holomorphic functions on $\OO$.
\end{example}

\section{The Kraft-Procesi Construction }
\label{sec:KP}

In this section we recall how  Kraft and Procesi in
\cite{KP1} and \cite{KP2} constructed the closures of complex classical
nilpotent orbits. We add to their construction  the  framework
of algebraic symplectic reduction.

Let $V$ be a complex vector space with $\bb$ a   bilinear form on $V$.
Let $G$ be the symmetry group of $\bb$. We consider only
the following three cases.
\begin{itemize}
\item[(i)]   $\bb$ is identically zero. Then $G=GL(V,\C)$.
\item[(ii)]    $\bb$ is   nondegenerate and symmetric. Then $G=O(V,\C)$.
\item[(iii)]   $\bb$ is   nondegenerate and symplectic. Then $G=Sp(V,\C)$.
\end{itemize}

Choose  a nilpotent orbit $\OO$ of $G$ and an element $\la$ in $\OO$.
Then  $\la$ corresponds (via the trace functional on $\End V$) to  $x\in\g$;
so   $x$ is in particular  a nilpotent  endomorphism of $V$.
We note that $\OO$ is connected, and so $\CO$ is irreducible,
except in the following situation.  If  $G=O(V,\C)$ where $\dim V$
is even and also the Jordan block size partition of $\OO$
is  very even (i.e., all parts are even and occur with even multiplicities),
then $\OO$ has two connected components and $\CO$ has
two irreducible components. See \cite{KP2}, \cite[Chapter 5]{CM}.

Let $V_d$  be the image of $x^d$. Then
\begin{equation}\label{eq:VVV}
V=V_0\supset V_1\supset\cdots \supset V_r\supset V_{r+1}=0
\end{equation}
where $r$ is the largest number such that $x^r\neq 0$.
We define a complex  vector space $L$   by
\begin{equation}\label{eq:L=II}
L=L(V_0,V_1)\oplus L(V_1,V_2)\oplus\cdots\oplus L(V_{r-1},V_{r})
\end{equation}
where $L(V_{d-1},V_{d})$ is obtained in the following way.
If $G=GL(V,\C)$ then
\begin{equation}\label{eq:L=}
L(V_{d-1},V_{d})=\Hom(V_d,V_{d-1})\oplus \Hom(V_{d-1},V_{d})
\end{equation}
If $G=O(V,\C)$ or $G=Sp(V,\C)$, then
\begin{equation}\label{eq:L=I}
L(V_{d-1},V_{d})=\Hom(V_d,V_{d-1})
\end{equation}

Next we construct a complex Lie group $S$ of the form
\begin{equation}\label{eq:S=}
S=S_{1}\times S_{2}\times\cdots\times S_{r}
\end{equation}
where $S_{d}$ is obtained in the following way.
To begin with, we put $\bb_0=\bb$ and  $S_0=G$.
If $G=GL(V,\C)$ then for each $d$ we put $\bb_d=0$ and $S_{d}=GL(V_d,\C)$.
If $G=O(V,\C)$ or $G=Sp(V,\C)$, then $V_d$  admits an intrinsic
 nondegenerate
complex bilinear form $\bb_d$ and we define
 $S_d$  to be  the symmetry group of $\bb_d$.
In more detail,  $\bb_d$ is  the bilinear form on $V_d$ defined by
$\bb_d(x^d(u),x^d(v))=\bb(u,x^d(v))$. It turns out   that $\bb_d$ is
  nondegenerate.
If  $\bb_{d-1}$ is orthogonal then $\bb_d$ is  symplectic
and we   put    $S_d=Sp(V_d,\C)$. If $\bb_{d-1}$ is  symplectic
  then $\bb_d$ is
orthogonal  and  we put  $S_d=O(V_d,\C)$.

Next  we construct  commuting  actions of  $G$ and $S$  on $L$.
If $G=GL(V,\C)$, we make  $G$ and $S$  act by
\begin{equation}
\begin{split}\label{eq:IISact}
(g,s_1,\dots,s_r)&\bullet(A_1,B_1,A_2,B_2,\dots,A_r,B_r)\\
& =(gA_1s_1\inv,  s_1B_1g\inv,  s_1A_2s_2\inv,  s_2B_2s_1\inv, \dots,
s_{r-1}A_rs_r\inv,  s_{r}B_rs_{r-1}\inv)
\end{split}
\end{equation}
where $A_d\in\Hom(V_d,V_{d-1})$ and $B_d\in\Hom(V_{d-1},V_d)$. If
$G=O(V,\C)$ or $G=Sp(V,\C)$, we make $G$ and $S$  act  by
\begin{equation}\label{eq:ISact}
(g,s_1,\dots,s_r)\bullet(C_1,C_2,\dots,C_r)=
(gC_1s_1\inv,   s_1C_2s_2\inv,   \dots, s_{r-1}C_rs_r\inv)
\end{equation}
where $C_d\in\Hom(V_d,V_{d-1})$.

$L$ has a  (complex) symplectic form $\Omega$ given by
$\Om=\Om_1+\Om_2+\cdots+\Om_r$ where $\Om_d$ is the symplectic
form on  $L(V_{d-1},V_{d})$ defined in the following way. If
$G=GL(V,\C)$, then $\Om_d(A+B,A'+B')=\tr(AB')-\tr(BA')$. If
$G=O(V,\C)$ or $G=Sp(V,\C)$, then $\Om_d(C,C')=\tr(C^*C)$ where
$C^*\in\Hom(V_{d-1},V_d)$ is the adjoint of $C$ defined by
$\bb_{d-1}(u,C^*(v))=\bb_d(C(u),v)$.

Now  $G$ and $S$ act faithfully and symplectically on $L$. Thereby
$G$ and $S$ identify with  commuting subgroups of the symplectic group
$Sp(L,\C)$. The action of  $Sp(L,\C)$ is Hamiltonian
with    canonical moment map
\begin{equation}\label{eq:mom}
L\to\spL^*
\end{equation}
Hence our actions of   $G$ and $S$ are
Hamiltonian  with induced moment maps (obtained by projection)
\begin{equation}\label{eq:ga_sig}
\ga:L\to\g^* \qquad\mbox{and}\qquad  \sig:L\to\s^*
\end{equation}
Then $\ga$ is $G$-equivariant and $S$-invariant, and
$\sig$ is $S$-equivariant and $G$-invariant.

Here are the explicit formulas for $\ga$ and $\sig$.
We may write these as $\g$-valued and $\s$-valued maps,
with the convention  that the trace functional on $\End L$
induces isomorphisms $\g\iso\g^*$ and $\s\iso\s^*$.
If $G=GL(V,\C)$, then $\ga$ and $\sig$  are given by
\begin{equation}\label{eq:II_mom}
(-A_1B_1) \qquad\mbox{and}\qquad
(B_1A_1-A_2B_2,\dots,   B_{r-1}A_{r-1}-A_rB_r,  B_rA_r)
 \end{equation}
If $G=O(V,\C)$ or $G=Sp(V,\C)$, then $\ga$ and $\sig$  are given by
\begin{equation}\label{eq:I_mom}
(-C_1C^*_1) \qquad\mbox{and}\qquad
(C_1^*C_1-C_2C^*_2,\dots, C_{r-1}^*C_{r-1}-C_rC_r^*,  C_r^*C_r)
 \end{equation}

\begin{remark}\label{rem:Howe}
For each $d=1,\dots,r$, $S_{d-1}$ and $S_d$ act on $L(V_{d-1},V_d)$ as
a  Howe dual pair.
\end{remark}

Let $\PP$ be the algebra of polynomial functions on $L$.
Then $\PP=S(L^*)$ is a graded Poisson algebra with respect to the
Poisson bracket  $\{\cdot,\cdot\}$ defined by $\Om$. Our grading is
$\PP=\oplus_{j\in\hN}\,\PP[j]$ where $\PP[j]=S^{2j}(L^*)$
and then $\{\PP[j],\PP[k]\}\subseteq\PP[j+k-1]$.
(This choice of halving the natural degrees is convenient as the aim is
to obtain $\RR$ as a subquotient of $\PP$.)
The momentum  functions $\ga_y$ ($y\in\g$) and $\sig_x$ ($x\in\s$)
are the component functions of $\ga$ and $\sig$; i.e.,
$\ga_y(m)=\langle\ga(m),y\rangle$  and
$\sig_x(m)=\langle\sig(m),x\rangle$ where $m\in L$.
The  $\ga_y$   and $\sig_x$  lie in $\PP[1]$ and satisfy the bracket relations
$\{\ga_y,\ga_{y'}\}=\ga_{[y,y']}$, $\{\sig_x,\sig_{x'}\}=\sig_{[x,x']}$ and
$\{\ga_y,\sig_x\}=0$.

Let $\II$ be the ideal in $\PP$ generated by the momentum  functions
$\sig_x$  where $x\in\s$.
Then $\II=\oplus_{j\in\hN}^\infty\,\II[j]$ is a graded ideal stable under
both $G$ and $S$.
Hence the quotient algebra $\PP/\II$ is a graded algebra on which
$G$ and $S$ act  by graded algebra automorphisms. The grading is
\[  \PP/\II=\oplus_{j\in\hN} (\PP/\II)^j   \]
where $(\PP/\II)^j=\PP[j]/\II[j]$.
Kraft and Procesi proved  that $\II$ is the full ideal of functions vanishing on
$\sig\inv(0)$. Thus $\PP/\II$ is the  coordinate ring   $\C[\sig\inv(0)]$
of the zero locus   of  $\sig$.

The \emph{algebraic symplectic reduction} $L\red$ of $L$
by $S$ is the Mumford quotient of  $\sig\inv(0)$ by $S$.
Thus $L\red$  is the affine complex algebraic variety with coordinate ring
\begin{equation}\label{eq:}
\C[L\red]=(\PP/\II)\INV=\PV/\IV
\end{equation}
where the superscript $inv$ denotes taking $S$-invariants.
Moreover  $\IV$ is a Poisson ideal in $\PV$.
So  $\C[L\red]$ inherits the structure of a  graded Poisson algebra.

Notice that $S$ contains the center $\Z_2=\{1,-1\}$ of $\SpL$
and the action of $\Z_2$ induces the decomposition
$\PP=\Pe\oplus\Po$ where the even part is the graded Poisson algebra
\begin{equation}\label{eq:Pe=}
\Pe=\oplus_{d\in\N}\,\PP[d]
\end{equation}
So  $\PV$ lies in   $\Pe$, and thus $\PV$ and $\C[L\red]$ are $\N$-graded.

The symplectic version of the Kraft-Procesi result is
\begin{theorem}\label{thm:KP}  \cite[Theorem 3.3]{KP1}, \cite[Theorem 5.3]{KP2}
The algebra homomorphism $\ga^*:\Sg\to\PP$ defined by
$y\mapsto\ga_y$ \textup{(}$y\in\g$\textup{)}
induces a $G$-equivariant  isomorphism of  $\N$-graded  Poisson algebras from
$\RR$ onto  $\PIV$.
\end{theorem}

The cited results of Kraft and Procesi are given in geometric language, and  the
reader who wants to read all the proofs in \cite{KP1} and
\cite{KP2} will need some knowledge
in algebraic geometry.  The statements in \cite[Theorem 3.3]{KP1} and
\cite[Theorem 5.3]{KP2} are easy to translate into algebra though,
since we are dealing with
\emph{affine} varieties.      Kraft and Procesi  show that
$\ga$ maps $\sig\inv(0)$  onto  $\CO$,  and moreover this surjection
 $\ga':\sig\inv(0)\to\CO$ is a quotient map for the action of $S$. In this setting
of   a reductive group acting on an affine variety, ``quotient map"  has a
very strong meaning coming from   Mumford's geometric invariant theory, as
explained in  \cite[\S1.4]{KP1} and \cite[\S0.11]{KP2}.
Precisely,   $\ga'$ being a quotient map means that the corresponding map
$\RR\to\PP/\II$ on coordinate rings  is injective and has image equal to
$(\PP/\II)\INV$.

In symplectic language then, Kraft and Procesi proved that
the moment map $\ga$ induces a $G$-isomorphism  of   affine
complex algebraic varieties  from $L\red$ onto $\CO$.
This isomorphism is also  equivariant with respect to the natural  $\C^*$-actions
on $L\red$ and $\CO$.
Finally,  since $\ga$ is a moment map it follows that   $\ga^*$
preserves the Poisson brackets. Thus we get Theorem \ref{thm:KP}.

\section{Weyl algebra $\W$ for  $L$}
\label{sec:Weyl}

The Kraft-Procesi construction realized $\RR$ as a subquotient,
namely   $\PIV$,  of  $\Pe$.
Our aim is to make a noncommutative analog of their
construction.

The image of the moment map (\ref{eq:mom}) is the closure $\oY$ of the
minimal nilpotent   orbit $\Y$ of $\SpL$, and  $\Pe=\C[\oY]$
as graded Poisson algebras.
We know by Example \ref{exam:min} that $\oY$ has a unique Dixmier
algebra $(\D,\xi_{\D},\tau_{\D},\vth_{\D})$, and then $\D$ is the quotient
of $\U(\spL)$ by its Joseph ideal.  In this section we will give a more
concrete model for this Dixmier algebra.  Then in \S\ref{sec:B}
we will perform  the noncommutative analog  of reduction.

Let $\W$ be the Weyl algebra for   $L^*$.
This means that $\W$ is the quotient of the tensor algebra  of $L^*$ by the
two-sided ideal generated by the elements
$a\otimes b-b\otimes  a-\PB{a}{b}$ where $a$ and $b$ lie in $L^*$.
Let $a\mapsto\hat{a}$ be the natural map $L^*\to\W$. We can identify
$\spL$ with $S^2L^*$ and then we have the Lie algebra embedding
\begin{equation}\label{eq:PP2->W}
\xi:\spL\longrightarrow  \W, \quad \xi^{ab}=\hat{a}\hat{b}+\hat{b}\hat{a}
\end{equation}

There is  an increasing algebra filtration  $\W=\cup_{j\in\hN}\;\W_j$
where $\W_j$ is the image of the space of tensors of degree at most $2j$.
We have $[\W_j,\W_k]\subset\W_{j+k-1}$.  Thus the associated graded algebra
$\gr\W=\oplus_{j\in\hN}\W_j/W_{j-\half}$ is  commutative and  the
commutator in $\W$ induces a Poisson bracket (of degree $-1$) on $\gr\W$.
In this way $\gr\W$ becomes a graded Poisson algebra.
Then  $\gr\W$ identifies naturally with $\PP$.

The symplectic group   $\SpL$ acts naturally on $\W$ by algebra
automorphisms. This action respects $\xi$, the filtration on $\W$,
the Poisson bracket on $\gr\W$, etc. The corresponding action of
$\spL$ on $\W$ is given by the operators $[\xi^{ab},\cdot]$.

We next choose a Cartan involution  $\vsig$ of $\spL$.
To do this, we go back into the Kraft-Procesi construction.
Recall that each space $V_d$ in (\ref{eq:VVV}) carried a
bilinear form $\bb_d$   ($d=0,\dots,r$).
We can choose a positive definite hermitian form $\bh_d$ on $V_d$
which is compatible with $\bb_d$ in the sense that the intersection of
$S_d$ with the unitary group of $\bh_d$ is a maximal compact subgroup
$K_d$ of $S_d$. (This is an equivalent version of the setup in \cite{KS}.)
These $\bh_d$ determine naturally a positive definite hermitian form
$\bh$ on $L$.  Now we define $\vsig(T)=-T^\dagger$  for $T\in\spL$,
where  $T^\dagger$  is the adjoint  of $T$ with respect to $\bh$.

Corresponding to $\vsig$ is a compact real form $\Spc$ of $\SpL$
with Lie algebra $\spc$.
For later use (see \S\ref{sec:B}), we notice that
$G\cap\Spc=K_0$ and $S\cap\Spc=K_1\times\cdots\times K_r$
are compact real forms of $G$ and $S$, which we will denote by
$\Gc$ and $\Sc$.

Now  $\W$ is a $(\spL\oplus\spL,\Spc)$-module, where the representation
\begin{equation}\label{eq:spsp}
\spL\oplus\spL\to\End\W
\end{equation}
is given by $(x,y)\cdot A=\xi^xA-A\xi^y$ and the action of $\Spc$ corresponds
to the subalgebra $\{(x,x): x\in\spc\}$.
The action  of  the center $\Z_2$ of $\Spc$   produces  the  decomposition
\begin{equation}\label{eq:W=}
\W=\We\oplus\Wo
\end{equation}
where $\We$ is  space of invariants for $\Z_2$.
The induced  filtration on $\We$ satisfies $\We_{p+\half}=\We_p$ if $p\in\N$.
So we might  as well  just consider the  algebra filtration
\begin{equation}\label{eq:We=}
 \We=\cup_{d\in\N}\,\We_d
\end{equation}
Now we  can make $\We$ into a Dixmier algebra.

\begin{proposition}\label{prop:We}
The  Dixmier algebra for the closure $\oY$ of the minimal nilpotent  orbit  of $\SpL$
is the quadruple $(\We,\xi,\tau,\vth)$,  for some unique choices of
$\tau$ and $\vth$.
\end{proposition}
\begin{proof}
The map $\xi$ induces  a filtered algebra isomorphism
$\pi:\U(\spL)/\JJ\to\We$ where $\JJ$ is   the Joseph ideal
(see \cite[\S5]{AB:starmin}).  Clearly  $\gr\We=\Pe$
and so $\pi$ induces a graded   isomorphism
$S(\spL)/\ICO[\oY]\iso\Pe$.
Now everything follows by Example \ref{exam:min}.
\end{proof}

\section{Dixmier Algebra for  $\CO$ }
\label{sec:B}
 $\We$ is, by means of (\ref{eq:spsp}), both a $(\gog,\Gc)$-module and
an $(\sos,\Sc)$-module,   and these two actions   commute.
\begin{definition}\label{def:B}
$\B$ is the $(\gog,\Gc)$-module obtained by taking the coinvariants
of $\We$ in the category of $(\sos,\Sc)$-modules.
\end{definition}

This means that $\B$ is the quotient $\We/\M$ where
$\M$ is the subspace spanned by  all
$\xi^xA-A\xi^y$  and  $A-s\cdot A$ where $x,y\in\s$, $s\in\Sc$ and $A\in\We$
(see \cite[Chapter II]{KV}).
Then $\B$ inherits from $\We$ an increasing $\Gc$-stable vector space
filtration   $\B=\cup_{d\in\N}\;\B_d$.

Let  $\WV$ be the algebra of  invariants for $\Sc$. Then $\WV$
lies in $\We$ (since $\Sc$ contains $\Z_2$) and so $\WV$
inherits from $\We$  an algebra filtration $\WV=\cup_{d\in\N}\,\WV_d$.
\begin{lemma}\label{lem:B}
The natural map
\begin{equation}\label{eq:WV->B}
\phi:\WV\to\B
\end{equation}
is surjective in each filtration degree and its kernel  is a two-sided ideal.
In this way,  $\B$ becomes  a filtered algebra.
The corresponding map $\gr\phi:\PV\to\gr\B$ is a surjective homomorphism of
graded  Poisson algebras.
\end{lemma}
\begin{proof}
We prove this in \S\ref{sec:lem}.
\end{proof}

Our main result is
\begin{theorem}\label{thm:B}
We have $\gr\B=\PV/\IV$. So $\gr\B\simeq \RR$ as graded Poisson algebras.
\end{theorem}
\begin{proof}
The proof occupies \S\ref{sec:thm}.
\end{proof}

\begin{corollary}\label{cor:B_Dix}
The quadruple $(\B,\xiB,\tauB,\vthB)$ is a Dixmier algebra for $\CO$, where
$\xiB$, $\tauB$ and $\vthB$ are the maps induced by $\xi$, $\tau$, and $\vth$.
\end{corollary}
\begin{proof}
We prove this in     \S\ref{sec:cor}.
\end{proof}

Let $J$ be the kernel of the algebra homomorphism $\wt{\xiB}:\Ug\to\B$
defined by $\xiB$.  Proposition \ref{prop:3.1} gives
\begin{corollary}\label{cor:B_J}
Suppose we exclude the cases where  $\OO$ is disconnected
\textup{(}so where $G=O(2n,\C)$ and  the Jordan block size partition
of $\OO$ is very even\textup{)}.
Then $J$ is a completely prime primitive ideal of $\Ug$ with
$\gr J=\ICO$  $\taug(J)=J$, and $\vthg(J)=J$.
\end{corollary}
The  methods we have used thus far give no information about the excluded cases.

\begin{remark}\label{rem:GL}
(i)  Suppose  $G=GL(n,\C)$. Then the space $L$ has a
$\GtS$-invariant polarization, and using this  we can describe $\B$ and $\xi$ in the
following way.  Let $X$ be the flag manifold of $G$ of flags
of the type in (\ref{eq:VVV}).  Let  $\D^{\half}(X)$ be the algebra of  twisted
differential operators for the (locally defined) square root of the canonical bundle
on $X$  as in \cite{Brylinski:EDQ}.
We can show (\cite{Brylinski:Dix_KP_II}) that $\B$ identifies with
$\D^{\half}(X)$ in such a way that  $\xi$ corresponds to the canonical mapping of
$\g$ into $\D^{\half}(X)$. Then  by \cite[Corollary 8.5]{Brylinski:EDQ},
$J$ is a maximal ideal in $\Ug$.

(ii) Suppose $G=O(n,\C)$ or $G=Sp(2n,\C)$. If $\OO$ is the minimal nilpotent
orbit, then $\B$ is the quotient of $\Ug$ by its Joseph ideal. This follows
by the result   in Example \ref{exam:min}.
\end{remark}

\section{Proof of Lemma \ref{lem:B}}
\label{sec:lem}

The action of $\Sc$ on $\W$ is completely reducible and locally finite, and
$S$ and $\Sc$ have the same invariants and the same irreducible subspaces.
So we can form the decomposition
\begin{equation}\label{eq:X}
\We=\WV\oplus\X
\end{equation}
where $\X$ is the sum of all   non-trivial $\Sc$-isotypic components.
Then $\X$ is the span of the elements $A-s\cdot A$ where
$A\in\We$ and $s\in\Sc$.  Then $\M=\MV\oplus\X$.
Hence the natural map  $\phi$  is surjective  and its kernel is $\MV$.
I.e., we have vector space isomorphisms
\begin{equation}\label{eq:isoB}
\WV/\MV\iso\We/\M\iso\B
\end{equation}
The decomposition (\ref{eq:X}) is compatible with the filtration on $\We$,
since the filtration is $\Sc$-invariant. Consequently $\phi$ is surjective in
each filtration degree.  So $\We$ and $\WV$ induce the same filtration on $\B$.

Next we show that $\MV$ is a two-sided ideal in $\WV$. To begin with,
$\MV$ lies inside the subspace $\M'$ of $\M$  spanned by all
$\xi^xA$ and $A\xi^x$ where $x\in\s$ and $A\in\We$.
This follows using (\ref{eq:X}).
So it suffices to show that $D\M'$ and $\M'D$ lie in $\M'$ if $D\in\WV$.
Obviously $DA\xi^x$ lies  in $\M'$.   Invariance of $D$ gives
$\xi^xD-D\xi^x=0$ and so $D\xi^xA=\xi^xDA$ lies  in $\M'$.
Thus $D\M'\subseteq\M'$;  similarly $\M' D\subseteq\M'$.

The associated graded algebra  $\gr\B$ is the  quotient  $\gr\WV/\gr\MV$.
We find $\gr\WV=\PV$. Now the  final assertion is clear.

\begin{remark}\label{rem:W_We}
If we replace $\We$ by $\W$ in Definition \ref{def:B},
then we get the same thing.
I.e., if    $\wt{\B}$   is  the  module of coinvariants of $\W$, then
$\wt{\B}$ identifies naturally with $\B$.
Indeed  $\wt{\B}=\W/\wt{\M}$ where $\wt{\M}$ is the subspace spanned by  all
$\xi^xA-A\xi^y$  and  $A-s\cdot A$ where now  $A\in\W$.
But then  $\wt{\M}=\Wo\oplus\M$ and so the natural map $\wt{\phi}:\WV\to\wt{\B}$
is surjective with the same kernel $\MV$.
Also the filtration $\wt{\B}=\cup_{j\in\hN}\, \wt{\B}_{j}$ induced by $\W$
reduces to the one induced by $\WV$ in the sense that
$\wt{\phi}(\WV_d)=\wt{\B}_d=\wt{\B}_{d+\half}$.
\end{remark}

\section{Proof of Theorem \ref{thm:B}}
\label{sec:thm}

We will compute $\gr\B$ by using a homology spectral sequence. We will
consider the relative Lie algebra homology $H(\sos,\Sc;\We)$. By definition
(see \cite[Chapter II, \S6-7]{KV}),
$H_j(\sos,\Sc;\We)$ is the $j$th derived functor,
in the category of $(\sos,\Sc)$-modules,  of the coinvariants.  So
\begin{equation}\label{eq:1}
\B=H_0(\sos,\Sc;\We)
\end{equation}
The idea is that we will introduce a filtration of the complex
that computes the homology  in such
a way that the induced filtration on  $H_0(\sos,\Sc;\We)$ is the one we
have already defined on $\B$. Then we will use the usual spectral sequence
of a filtered complex to compute $\gr H_0(\sos,\Sc;\We)$.
The computation  will rely on the geometric   result of Kraft and Procesi
that (in the notation of  \S\ref{sec:KP}) $\sig\inv(0)$ is a complete intersection.

To begin with, we have  $\sos=\kk\oplus\p$ where
$\kk=\{(x,x): x\in\s\}$ and    $\p=\{(x,-x): x\in\s\}$.
Then $\kk$ is the complexified Lie algebra of $\Sc$.  The standard complex
(\cite[page 163]{KV}) for computing $H(\sos,\Sc;\We)$  is
 \begin{equation}\label{eq:2}
0\mapleft{} \wedge^{0}\p\otimes_{\Sc}\We
\mapleft{\del} \wedge^{1}\p\otimes_{\Sc}\We
\mapleft{\del} \cdots
\mapleft{\del} \wedge^{m}\p\otimes_{\Sc}\We\mapleft{} 0
\end{equation}
Here $m=\dim\s$ and  $\otimes_{\Sc}$ denotes the $\Sc$-coinvariants
of  the tensor product.       We call this complex $A$ where
$\leftsuper{t}A=\wedge^{t}\p\otimes_{\Sc}\We$.

The differential  $\del$ in (\ref{eq:2})  is given  by
\begin{equation}\label{eq:3}
\del(Y_1\wedge\cdots\wedge Y_t\otimes D)=
\sum_{l=1}^t(-1)^lY_1\wedge\cdots\wedge\widehat{Y_l}
\wedge\cdots\wedge Y_t
\otimes \Pi^{Y_l}(D)
\end{equation}
where the $Y_i$ lie in $\p$,  $D\in\We$, and $\Pi$ is the representation
(\ref{eq:spsp}).       (Notice the  terms involving $[Y_i,Y_j]$
are not present  because $[\p,\p]\subseteq\kk$.)
To make the complex more transparent, we identify $\p$ with $\s$
so that  $(x,-x)$ corresponds to $x$. Then (\ref{eq:2}) becomes
\begin{equation}\label{eq:4}
\del(x_1\wedge\cdots\wedge x_t\otimes D)=
\sum_{l=1}^t(-1)^lx_1\wedge\cdots\wedge\widehat{x_l}
\wedge\cdots\wedge x_t
\otimes (\xi^{x_l}D+D\xi^{x_l})
\end{equation}

Next we  define an increasing filtration of $A$ by the spaces
\begin{equation}\label{eq:5}
\leftsuper{t}A^d=\wedge^{t}\p\otimes_{\Sc}\We_{d-t}
\end{equation}
where we set $\We_j=0$ if $j<0$.
Then   $\del\,\leftsuper{t}A^d\subseteq \leftsuper{t-1}A^d$. This follows since
the $\xi^x$ lie in $\We_1$ and so $\p\cdot\We_d\subseteq\We_{d+1}$.
So we have in hand a  filtration of the complex (\ref{eq:1}).
We put $A^{d,q}=\leftsuper{d+q}A^d$; then $d$ is the filtration degree and
$q$ is the complementary degree. The induced   filtration  on the homology is
$H(\sos,\Sc;\We)=\cup_{d\in\N}\,F^d$
where $F^d=\oplus_{q\in\Z}\, F^{d,q}$ and
$F^{d,q}$ is the $d$th filtration piece of  $H_{d+q}(\sos,\Sc;\We)$. The
associated graded space $\gr\,H(\sos,\Sc;\We)$
is the direct sum of the spaces
\begin{equation}\label{eq:6}
\gr^d\,H_{d+q}(\sos,\Sc;\We)=F^{d,q}/F^{d-1,q+1}
\end{equation}

Notice that $\leftsuper{0}A=\We$ and the  filtration on $\leftsuper{0}A$
defined by (\ref{eq:5}) is  the same one  as in (\ref{eq:We=}).
So $\gr\,H_{0}(\sos,\Sc;\We)=\gr\,\B$.
Our goal is to prove
\begin{equation}\label{eq:7}
\gr^d\,H_{0}(\sos,\Sc;\We)=(\PIV)^d
\end{equation}

Now we consider  the spectral sequence $E_0,E_1,\dots$
associated to our  filtered complex.
(See e.g., \cite[Appendix D]{KV} or \cite[Chapter I,\S4]{Godement}   for
the construction of this spectral sequence in  the general setting.)
The $E_0$  term is given by  $E_0^{d,q}=A^{d,q}/A^{d-1,q+1}$  and so
\begin{equation}\label{eq:8}
E_0^{d,q}=\wedge^{d+q}\,\p\,\otimes_{\Sc}\,\We_{-q}/\We_{-q-1}
\end{equation}
The identification  $\gr\We=\Pe$ gives
\begin{equation}\label{eq:9}
E_0^{d,q}=\wedge^{d+q}\,\p\,\otimes_{\Sc}\,\PP[-q]
\end{equation}
(Thus the $E_0$ term occupies the octant of  the $d,q$
plane where $q\le 0$ and $d+q\ge 0$.)
So   $E^d_0$ is the complex
\begin{equation}\label{eq:10}
0\mapleft{} \wedge^{0}\p\otimes_{\Sc}\PP[d]
\mapleft{\del_0} \wedge^{1}\p\otimes_{\Sc}\PP[d-1]
\mapleft{\del_0}\cdots \mapleft{\del_0}
\wedge^{m}\p\otimes_{\Sc}\PP[d-m]\mapleft{}0
\end{equation}
The boundary $\del_0$ is induced by $\del$.
We can easily compute $\del_0$ since the  natural projection maps
$\psi_d:\W_d\to\PP[d]$   are given by
$\psi_d(\hat{a}_1\cdots\hat{a}_{2d})=a_1\cdots a_{2d}$
where $a_i\in L^*$ (cf. \S\ref{sec:Weyl}).
So for  $D\in\W_d$ we have
$\psi_{d+1}(\xi^xD)=\psi_{d+1}(D\xi^x)=\sig_x\psi_d(D)$.
Thus  (\ref{eq:4}) gives
\begin{equation}\label{eq:11}
\del_0(x_1\wedge\cdots\wedge x_t\otimes f)=
\sum_{l}(-1)^lx_1\wedge\cdots\wedge\widehat{x_l}\wedge\cdots\wedge x_t
\otimes (2\sig_{x_l}f)
\end{equation}

The total complex $E_0$ is
\begin{equation}\label{eq:12}
0\mapleft{} \wedge^{0}\p\otimes_{\Sc}\,\Pe
\mapleft{\del_0} \wedge^{1}\p\otimes_{\Sc}\,\Pe\mapleft{\del_0}
\cdots \mapleft{\del_0}\wedge^{m}\p\otimes_{\Sc}\,\Pe\mapleft{} 0
\end{equation}
The homology $H(E_0)$, together with a  differential $\del_1$,
is the $E_1$ term of the spectral sequence.
More precisely,   $E_1^{d,q}=H_{d+q}(E_0^{d,*})$.

To compute $E_1$, we observe that $H(E_0)$ is  the  $\Sc$-coinvariants
of the homology  of  the complex
\begin{equation}\label{eq:13}
0\mapleft{} \wedge^{0}\p\otimes\,\Pe
\mapleft{\del_0} \wedge^{1}\p\otimes\,\Pe\mapleft{\del_0}
\cdots \mapleft{\del_0}\wedge^{m}\p\otimes\,\Pe\mapleft{} 0
\end{equation}
Indeed,  $E_0$ is the $\Sc$-coinvariants of (\ref{eq:13}), and
taking coinvariants commutes with taking homology.  The latter
follows because each space $\wedge^{t}\p\otimes\,\Pe$ is a locally finite
$\Sc$-representation, and for any such representation $\V$,  the natural map
$\V^{\Sc}\to\V_{\Sc}$  from invariants  to coinvariants is an isomorphism.

To compute the homology of (\ref{eq:13}), we recognize (\ref{eq:13})  as the
Koszul complex $K$  of the sequence $\sig_{y_1},\dots,\sig_{y_m}$ in $\Pe$
where $y_1,\dots,y_m$ is any  basis of $\s$.  Recall from
\S\ref{sec:KP} that $\II$ is the ideal in $\PP$ generated by the $\sig_{y_i}$.
Kraft and Procesi proved in \cite[Theorem 3.3]{KP1} and \cite[Theorem 5.3]{KP2}
that  the subscheme $\sig\inv(0)$ of $L$ is a reduced complete intersection,
i.e.,  $\sig_{y_1},\dots,\sig_{y_m}$  is  a regular sequence in $\PP$.
Let us consider the Koszul complex $\wt{K}$ of this sequence   in $\PP$.
By a well known result of commutative algebra
(see  \cite[III, Proposition 7.10A]{Hartshorne})
the homology of  $\wt{K}$   is concentrated in degree zero and
$H_0(\wt{K})=\PP/\II$  as graded algebras.
But $K$ is simply obtained from $\wt{K}$  by taking $\Z_2$-invariants.
Hence the homology of  $K$   is concentrated in degree zero and
$H_0(K)=\Pe/\Ie$  as graded algebras.
Then the module  $H_0(K)_{\Sc}$ of coinvariants  identifies with $\PIV$.

Thus $E_1$ is the complex
\begin{equation}\label{eq:14}
0\mapleft{}\PIV\mapleft{\del_1}0\mapleft{\del_1}\cdots  \mapleft{\del_1}0
\mapleft{}0
\end{equation}
where the differentials $\del_1^{d,q}:E_1^{d,q}\to E_1^{d-1,q}$ are obviously zero and
\begin{equation}\label{eq:15}
 E_1^{d,-d}=(\PIV)^d \qquad\mathrm{while}\qquad E_1^{d,q}=0\;\mathrm{if}\;  q\neq -d
\end{equation}

Now we can compute the rest of the spectral sequence.  We know
$E_{r+1}$ is the homology of $E_{r}$ with respect to a differential  $\del_r$;
i.e. $E_{r+1}^{d,q}=\mathrm{ker}\,\del_{r}^{d,q}/\mathrm{im}\,\del_{r}^{d+r,q-r+1}$
where $\del_r^{d,q}$ maps $E_r^{d,q}$ to $E_r^{d-r,q+r-1}$.
For   $r\ge 1$, we find   that   $E_1^{d,q}=E_r^{d,q}$
and the differentials  $\del_r^{d,q}$ are  all zero.

The $\Ei$ term of the spectral sequence  satisfies
\begin{equation}\label{eq:16}
\Ei^{d,q}=\gr^d\,H_{d+q}(\sos,\Sc;\We)
\end{equation}
Our final step  is  to show  our spectral sequence converges in  that
\begin{equation}\label{eq:17}
E_1^{d,q}=\Ei^{d,q}
\end{equation}
This will finish off  the proof of Theorem \ref{thm:B}
because then (\ref{eq:15})    gives  the desired result (\ref{eq:7}).

The convergence   (\ref{eq:17})  follows formally from the two properties:
(i) $A^d=0$ if $d<0$ where $A^d=\oplus_{q\in\Z}A^{d,q}$
and (ii) $A^{d,q}$ has  finite dimension.
Indeed, following the notation in \cite[I,\S4.2]{Godement}, we have
\begin{equation}\label{eq:18}
\begin{array}{rcl}
E_r^{d,q}&=&Z_r^{d,q}/(B_{r-1}^{d,q}+Z_{r-1}^{d-1,q+1}) \\[6pt]
\Ei^{d,q}&=&\Zi^{d,q}/(\Bi^{d,q}+\Zi^{d-1,q+1})
\end{array}
\end{equation}
where   $Z_r^{d,q}=\{z\in A^{d,q} : \del z\in A^{d-r}\}$,
$\Zi^{d,q}=A^{d,q}\cap \ker\del$,
$B_r^{d,q}=A^{d,q}\cap\,\del A^{d+r}$ and $\Bi^{d,q}=A^{d,q}\cap\,\del A$.
Suppose we fix $d$ and $q$. Then (i) gives $Z_r^{d,q}=\Zi^{d,q}$
and $Z_{r-1}^{d-1,q+1}=\Zi^{d-1,q+1}$  if $r>d$,
and (ii) gives $B_{r-1}^{d,q}=\Bi^{d,q}$ if $r$ is large enough.
Therefore    $E_r^{d,q}=\Ei^{d,q}$    for $r$  large enough.
But we found $E_1^{d,q}=E_r^{d,q}$ ($r\ge 1$) and so  (\ref{eq:17}) follows.

We remark that it also follows that  $\gr H_j(\sos,\Sc;\We)=0$ for $j>0$.
Thus we have proven
\begin{proposition}\label{prop:H_j}
$H_j(\sos,\Sc;\We)=0$ if $j>0$.
\end{proposition}

\section{Proof of Corollary \ref{cor:B_Dix}}
\label{sec:cor}

$\Sc$ and $\Gc$ are commuting subgroups of $Sp(L)$.
It follows that $\xi$ maps $\g$ into $\WV$ and so $\xi$ induces $\xiB$
where $\xiB$ is the composition $\g\mapright{\xi}\WV\to\B$.
Next $\tau$ is $Sp(L)$-invariant and so in particular  is $\Sc$-invariant.
Then $\tau$ preserves $\WV$. We see that $\tau$ preserves $\M$, and
so $\tau$ preserves $\MV$. Hence
$\tau$ induces a filtered algebra anti-involution $\tauB$ of $\B=\WV/\MV$.
Finally, $\vth$  is $\Spc$-invariant and so in particular  is $\Sc$-invariant.
Then $\vth$ preserves $\WV$.   We see that $\vth$ preserves $\M$, and
so $\vth$ induces an antilinear  filtered algebra involution $\vthB$ of $\B$.

Thus we have in place our Dixmier algebra data for $\CO$. It is clear
because of Proposition \ref{prop:We}   that the axioms are satisfied.

\bibliographystyle{plain}

\end{document}